\documentclass[12pt,a4paper]{amsart}

\usepackage{t1enc}
\usepackage[latin1]{inputenc}
\usepackage[english]{babel}
\usepackage{hyperref}
\usepackage{graphicx}

\usepackage{latexsym}
\usepackage{amsmath,amsfonts,amssymb,amsthm}

\begin{document}

\title{On Cimmino Integrals as Residues of Zeta Functions}

\author{Sergio Venturini}

\address{
	S. Venturini:
	Dipartimento Di Matematica,
	Universit\`{a} di Bologna,
	\,\,Piazza di Porta S. Donato 5 ---I-40127 Bologna,
	Italy}
\email{venturin@dm.unibo.it}



%
%
\def\R{{\rm I\kern-.185em R}}
\def\RR{\mathbb{R}}
\def\C{{\rm\kern.37em\vrule height1.4ex width.05em depth-.011em\kern-.37em C}}
\def\CC{\mathbb{C}}
\def\N{{\rm I\kern-.185em N}}
\def\NN{\mathbb{N}}
\def\Z{{\bf Z}}
\def\ZZ{\mathbb{N}}
\def\Q{{\mathchoice{\hbox{\rm\kern.37em\vrule height1.4ex width.05em 
depth-.011em\kern-.37em Q}}{\hbox{\rm\kern.37em\vrule height1.4ex width.05em 
depth-.011em\kern-.37em Q}}{\hbox{\sevenrm\kern.37em\vrule height1.3ex 
width.05em depth-.02em\kern-.3em Q}}{\hbox{\sevenrm\kern.37em\vrule height1.3ex
width.05em depth-.02em\kern-.3em Q}}}}
\def\P{{\rm I\kern-.185em P}}
\def\H{{\rm I\kern-.185em H}}
%
\def\Aleph{\aleph_0}
\def\ALEPH#1{\aleph_{#1}}
\def\sset{\subset}\def\ssset{\sset\sset}
%
\def\bar#1{\overline{#1}}
\def\dim{\mathop{\rm dim}\nolimits}
\def\half{\textstyle{1\over2}}
\def\Half{\displaystyle{1\over2}}
\def\mlog{\mathop{\half\log}\nolimits}
\def\Mlog{\mathop{\Half\log}\nolimits}
\def\Det{\mathop{\rm Det}\nolimits}
\def\Hol{\mathop{\rm Hol}\nolimits}
\def\Aut{\mathop{\rm Aut}\nolimits}
\def\Re{\mathop{\rm Re}\nolimits}
\def\Im{\mathop{\rm Im}\nolimits}
\def\Ker{\mathop{\rm Ker}\nolimits}
\def\Fix{\mathop{\rm Fix}\nolimits}
\def\Exp{\mathop{\rm Exp}\nolimits}
\def\sp{\mathop{\rm sp}\nolimits}
\def\id{\mathop{\rm id}\nolimits}
\def\Rank{\mathop{\rm rk}\nolimits}
\def\Trace{\mathop{\rm Tr}\nolimits}
\def\Res{\mathop{\rm Res}\limits}
\def\cancel#1#2{\ooalign{$\hfil#1/\hfil$\crcr$#1#2$}}
\def\prevoid{\mathrel{\scriptstyle\bigcirc}}
\def\void{\mathord{\mathpalette\cancel{\mathrel{\scriptstyle\bigcirc}}}}
\def\n{{}|{}\!{}|{}\!{}|{}}
\def\abs#1{\left|#1\right|}
\def\norm#1{\left|\!\left|#1\right|\!\right|}
\def\nnorm#1{\left|\!\left|\!\left|#1\right|\!\right|\!\right|}
%
\def\upperint{\int^{{\displaystyle{}^*}}}
\def\lowerint{\int_{{\displaystyle{}_*}}}
\def\Upperint#1#2{\int_{#1}^{{\displaystyle{}^*}#2}}
\def\Lowerint#1#2{\int_{{\displaystyle{}_*}#1}^{#2}}
%
\def\rem #1::#2\par{\medbreak\noindent{\bf #1}\ #2\medbreak}
\def\proclaim #1::#2\par{\removelastskip\medskip\goodbreak{\bf#1:}
\ {\sl#2}\medskip\goodbreak}
\def\ass#1{{\rm(\rmnum#1)}}
\def\assertion #1:{\Acapo\llap{$(\rmnum#1)$}$\,$}
\def\Assertion #1:{\Acapo\llap{(#1)$\,$}}
\def\acapo{\hfill\break\noindent}
\def\Acapo{\hfill\break\indent}
\def\proof{\removelastskip\par\medskip\goodbreak\noindent{\it Proof.\/\ }}
\def\prova{\removelastskip\par\medskip\goodbreak
\noindent{\it Dimostrazione.\/\ }}
\def\risoluzione{\removelastskip\par\medskip\goodbreak
\noindent{\it Risoluzione.\/\ }}
\def\qed{{$\Box$}\par\smallskip}
\def\BeginItalic#1{\removelastskip\par\medskip\goodbreak
\noindent{\it #1.\/\ }}
\def\iff{if, and only if,\ }
\def\sse{se, e solo se,\ }
\def\rmnum#1{\romannumeral#1{}}
\def\Rmnum#1{\uppercase\expandafter{\romannumeral#1}{}}
\def\smallfrac#1/#2{\leavevmode\kern.1em
\raise.5ex\hbox{\the\scriptfont0 #1}\kern-.1em
/\kern-.15em\lower.25ex\hbox{\the\scriptfont0 #2}}
%
\def\Left#1{\left#1\left.}
\def\Right#1{\right.^{\llap{\sevenrm
\phantom{*}}}_{\llap{\sevenrm\phantom{*}}}\right#1}
%
%
%
\def\dimens{3em}
\def\symb[#1]{\noindent\rlap{[#1]}\hbox to \dimens{}\hangindent=\dimens}
\def\references{\bigskip\noindent{\bf References.}\bigskip}
\def\art #1 : #2 ; #3 ; #4 ; #5 ; #6. \par{#1, 
{\sl#2}, #3, {\bf#4}, (#5), #6.\par\smallskip}
\def\book #1 : #2 ; #3 ; #4. \par{#1, {\bf#2}, #3, #4.\par\smallskip}
\def\freeart #1 : #2 ; #3. \par{#1, {\sl#2}, #3.\par\smallskip}
%
%
%
%
\def\name{\hbox{Sergio Venturini}}
\def\snsaddress{\indent
\vbox{\bigskip\bigskip\bigskip
\name
\hbox{Scuola Normale Superiore}
\hbox{Piazza dei Cavalieri, 7}
\hbox{56126 Pisa (ITALY)}
\hbox{FAX 050/563513}}}
\def\cassinoaddress{\indent
\vbox{\bigskip\bigskip\bigskip
\name
\hbox{Universit\`a di Cassino}
\hbox{via Zamosch 43}
\hbox{03043 Cassino (FR)}
\hbox{ITALY}}}
\def\bolognaaddress{\indent
\vbox{\bigskip\bigskip\bigskip
\name
\hbox{Dipartimento di Matematica}
\hbox{Universit\`a di Bologna}
\hbox{Piazza di Porta S. Donato 5}
\hbox{40127 Bologna (BO)}
\hbox{ITALY}
\hbox{venturin@dm.unibo.it}
}}
\def\homeaddress{\indent
\vbox{\bigskip\bigskip\bigskip
\name
\hbox{via Garibaldi, 7}
\hbox{56124 Pisa (ITALY)}}}
\def\doubleaddress{
\vbox{
\hbox{\name}
\hbox{Universit\`a di Cassino}
\hbox{via Zamosch 43}
\hbox{03043 Cassino (FR)}
\hbox{ITALY}
\smallskip
\hbox{and}
\smallskip
\hbox{Scuola Normale Superiore}
\hbox{Piazza dei Cavalieri, 7}
\hbox{56126 Pisa (ITALY)}
\hbox{FAX 050/563513}}}
\def\sergio{{\rm\bigskip
\centerline{Sergio Venturini}
\leftline{\bolognaaddress}
\bigskip}}
%
%
%
%
%
%
\def\a{\alpha}
\def\bg{\beta}
\def\g{\gamma}
\def\G{\Gamma}
\def\dg{\delta}
\def\D{\Delta}
\def\e{\varepsilon}
\def\eps{\epsilon}
\def\z{\zeta}
\def\th{\theta}
\def\T{\Theta}
\def\k{\kappa}
\def\lg{\lambda}
\def\Lg{\Lambda}
\def\m{\mu}
\def\n{\nu}
\def\r{\rho}
\def\s{\sigma}
\def\Sg{\Sigma}
\def\ph{\varphi}
\def\Ph{\Phi}
\def\x{\xi}
\def\om{\omega}
\def\Om{\Omega}


\newtheorem{theorem}{Theorem}[section]
\newtheorem{proposition}{Proposition}[section]
\newtheorem{lemma}{Lemma}[section]
\newtheorem{corollary}{Corollary}[section]
\newtheorem{remark}{Remark}[section]
\newtheorem{definition}{Definition}[section]

\newtheorem{teorema}{Teorema}[section]
\newtheorem{proposizione}{Proposizione}[section]
\newtheorem{corollario}{Corollario}[section]
\newtheorem{osservazione}{Osservazione}[section]
\newtheorem{definizione}{Definizione}[section]
\newtheorem{esempio}{Esempio}[section]
\newtheorem{esercizio}{Esercizio}[section]
\newtheorem{congettura}{Congettura}[section]

\bibliographystyle{abbrv}
\def\cspace#1{C(\mathbb{R}^{#1})}
\def\lspace#1{L^1(\mathbb{R}^{#1})}
\def\sspace#1{\mathcal{S}(\mathbb{R}^{#1})}
\def\cbsspace#1#2#3{\mathcal{S}_{#3+#1}^{#3+#2}(\mathbb{R}^{#3})}
\def\expa{\alpha}
\def\expb{\beta}

\def\ndim{n}
\def\rvar{x}
\def\fvar{y}
\def\cvar{s}
\def\tvar{t}
\def\ivar{\om}

\def\Lattice{L}
\def\lA{A}
\def\qA{Q}
\def\qB{S}
\def\qC{C}
\def\qc{c}
\def\gbell{g}
\def\gbella#1{\gbell_{#1}}
\def\gbellab#1#2{\gbell_{{#1},{#2}}}
\def\lgbellab#1#2#3{\gbell_{{#1},{#2},{#3}}}
\def\cbt#1#2#3{\theta_{#1}\left({#2},{#3}\right)}
\def\cbts#1#2#3{\theta_{#1}^*\left({#2},{#3}\right)}
\def\cbx#1#2#3{\xi_{#1}\left({#2},{#3}\right)}
\def\cbza#1#2{\zeta\left(\qform{#1},{#2}\right)}
\def\cbzab#1#2#3{\zeta\left(\qform{#1},\qform{#2},{#3}\right)}
\def\cblza#1#2#3{\zeta_{#1}\left(\qform{#2},{#3}\right)}
\def\cblzab#1#2#3#4{\zeta_{#1}\left(\qform{#2},\qform{#3},{#4}\right)}
\def\cbcza#1#2{\zeta\left({#1},{#2}\right)}
\def\cbczab#1#2#3{\zeta\left({#1},{#2},{#3}\right)}
\def\cbtsum#1#2{\sum_{{#1}\in\Z^{#2}}}
\def\cbtsumnz#1#2{{\sum_{{#1}\in\Z^{#2}}}^{'}}

\def\cblsum#1#2{\sum_{{#1}\in{#2}}}
\def\cblsumnz#1#2{{\sum_{{#1}\in{#2}}}^{'}}

\def\cbGff{G_{f}}
\def\cbZff{Z_{f}}
\def\cbGhf{G_{\hat{f}}}
\def\cbZhf{Z_{\hat{f}}}

\def\qform#1{q_{#1}}
\def\qformb#1#2{q_{{#1},{#2}}}
\def\eucdot[#1,#2]{
	\left\langle{#1},{#2}\right\rangle
}
\def\cbdual#1{{#1}^\wedge}
\def\cbmdual#1{\hat{#1}}

\def\cbTranspose#1{{#1}^t}
\def\negone{{{-1}}}

\begin{abstract}
The following paper is a variation on a theme of
Gianfranco Cimmino on some integral representation
formulas for the solution of a linear equations system.

Cimmino was probably motivated for giving a representation
formula suitable not only for theoretical investigations
but also for applied computation.

In this paper we will prove that the Cimmino integrals
are strictly related to the residues of some zeta-like
functions associated to the linear system.

\end{abstract}

\maketitle

\section{\label{section:Introduction}Introduction}
Gianfranco Cimmino was born in Naples on March 12, 1908.

He received his Laurea degree in Mathematics at the University of Naples 
under the direction of Mauro Picone (1885-1977) in 1927.

At the end of 1939 Cimmino moved permanently to the
University of Bologna to occupy the chair of Mathematical Analysis.

Cimmino died in Bologna on May 30, 1989.

Gianfranco Cimmino (1908-1989) was a student
of Renato Caccioppoli (1904-1959)
and Mauro Picone
toghether with Giuseppe Scorza Dragoni (1907-1996)
and Carlo Miranda (1912-1982).

Mauro Picone founded the
``Istituto per le Applicazioni del Calcolo'' (IAC) in 1927.
Indeed, long before the introduction of digital computers,
Picone had the intuition of the potential impact on real-life problems
of the combination of computational methods with mathematical abstraction.

Probably the influence of Picone ideas 
and Cesari papers
\cite{article:CesaruNumerA} and \cite{article:CesaruNumerB}
on numerical solution of linear systems leads
Cimmino to be interested to the numerical treatment of the solutions
of linear systems of algebraic equations.

Among his researches he also obtained an interesting representation
of the solution of a linear system of real linear equations
that now we describe.

Let
\begin{eqnarray}
\label{eq::baseLinearSystem}
\lA x&=&b
\end{eqnarray}
be a linear system of $\ndim$ equation and $\ndim$ unknown,
where the unknown values $x_1,\ldots,x_\ndim$ are
the components if the column vector $x\in\R^\ndim$
and $\lA$ is non singular matrix with real coefficient of
order $\ndim$.

The well known Cramer's rule say that if $D=\det(\lA)\neq0$ then
\begin{eqnarray}
\label{eq::CramerSolution}
x_i&=&\frac{D_i}{D},\quad i=1,\ldots\ndim,
\end{eqnarray}
where $D_i$ is the determinant of the matrix obtained replacing
the $i-$th column of the matrix $A$ with the column vector $b$.

The alternative representation of (\ref{eq::CramerSolution}) given by Cimmino is
\begin{eqnarray}
\label{eq::CimminoSolution}
x_i&=&\frac{C_i}{C},\quad i=1,\ldots\ndim,
\end{eqnarray}
where
\begin{eqnarray}
\label{eq::CimminoNum}
	C&=&
		\int_{S^{\ndim-1}}\norm{\cbTranspose\lA u}^{-\ndim}\,du,
\end{eqnarray}
and
\begin{eqnarray}
\label{eq::CimminoDen}
	C_i&=&
		\ndim\int_{S^{\ndim-1}}
			\norm{\cbTranspose\lA u}^{-\ndim-2}
				\eucdot[b,u]\eucdot[\cbTranspose\lA x,e_i]\,du.
\end{eqnarray}

The integration here is made with respect to the standard
$(\ndim-1)-$dimensional measure on $S^{\ndim-1}$, the boundary of the Euclidean
unit ball in $\R^\ndim$,
$\norm{\cdot}$ and $\eucdot[\cdot,\cdot]$ stand respectively for the
Euclidean norm and the Euclidean inner product on $\R^\ndim$
and $e_1,\ldots,e_\ndim$ is the canonical basis of $\R^\ndim$,
and $\cbTranspose\lA$ denote the transpose of the matrix $\lA$.

In \cite{article:CimminoA}, \cite{article:CimminoB} and
\cite{article:CimminoC} Cimmino gives a probabilistic argument which
justify the existence of such kind of formulas and in
\cite{article:CimminoD} he also give an elementary but not trivial
proof of (\ref{eq::CimminoSolution}).
See also (\cite{article:CimminoE}) for some applications.
We refer to \cite{article:BenziPerCimmino}
for a detailed discussion and a background information
on Cimmino's papers in the field of the numerical analysis.

The purpose of this paper is to show that the Cimmino ideas
fit nicely into the theory of the residues of zeta-like functions,
a powerful tool coming from the analytic number theory.

Namely we will show that
(\ref{eq::CimminoNum}) and (\ref{eq::CimminoDen})
are the integral representation of the residues of suitable
zeta-like function associate to the matrix $A$ and the vector $b$
of the linear system (\ref{eq::baseLinearSystem});
see Theorem \ref{stm::ResSolve} for the complete statement.

It should be noted that formulas \ref{eq::CimminoNum}
and \ref{eq::CimminoDen} actually are due to Jacobi: cf. \cite{article:JacobiIntegDet}.

\section{\label{section:ThetaFunctions}Theta functions and their Mellin transform}
Let us begin with the following (almost trivial) observation:

\begin{proposition}\label{stm::resA}
Let $f$ be a continuous complex function defined on the
real interval $[0,1]$.
Assume that for some constants $R,\alpha,\beta\in\R$ with $\alpha>\beta$
we have
\begin{eqnarray*}
	f(t)&=&R t^{-\alpha}+O(t^{-\beta}),\quad t\to0^+.
\end{eqnarray*}
Then, given $s\in\C$, the integral
\begin{eqnarray*}
	g(s)&=&\int_0^1f(t)t^s\,\frac{dt}{t}
\end{eqnarray*}
converges absolutely on the half space $\Re(s)>\alpha$
and extends to a meromorphic function on the half space $\Re(s)>\beta$
having a simple pole at $s=\alpha$ with residue $R$.
\end{proposition}

\proof
Since $f(t)=O(t^{-\alpha})$ then $g(s)$ converges and is holomorphic
when $\Re(s)>\alpha$. If we denote
\begin{eqnarray*}
	h(t)&=&f(t)-R t^{-\alpha}
\end{eqnarray*}
then 
\begin{eqnarray*}
	g(s)&=&\int_0^1f(t)t^s\,\frac{dt}{t}
	=\frac{R}{s-\alpha}+\int_0^1h(t)t^s\,\frac{dt}{t}.
\end{eqnarray*}
Since $h(t)=O(t^{-\beta})$ the last integral define a holomorphic
function on the half space $\Re(s)>\beta$ and we are done.

\qed

%

Thus the singularities of the analytic function $g(s)$
describe the behaviour of the function $f(t)$ as $t\to0^+$.

Let us recall that $\cspace\ndim$ and $\lspace\ndim$
denotes respectively the space of the continuous complex function
on $\R^\ndim$ and
the space of the absolutely integrable
complex functions with respect to the Lebesgue measure on $\R^\ndim$.

Given $f\in\lspace\ndim$ the \emph{Fourier transform} of $f$ is
the function $\hat{f}$ defined by the formula
\begin{eqnarray*}
	\hat{f}(y)&=&\int_{\mathbb{R}^\ndim}f(x)e^{-2\pi i<x,y>}\,dx.
\end{eqnarray*}
Observe that
\begin{eqnarray*}
	\hat{f}(0)&=&\int_{\mathbb{R}^\ndim}f(x)\,dx.
\end{eqnarray*}

We denote by $\sspace\ndim$ is the Schwartz space of smooth functions
rapidly decreasing at the infinity together with their derivatives
of all orders.
Given $\expa,\expb>0$ two positive real number we denote by
\begin{eqnarray*}
	&&\cbsspace\expa\expb\ndim
\end{eqnarray*}
the space of all measurable function such that
\begin{eqnarray*}
	&&\sup_{x\in\R^\ndim}\abs{f(x)}(1+\norm{x}^{\ndim+\expa})<+\infty
\end{eqnarray*}
and
\begin{eqnarray*}
	&&\sup_{y\in\R^\ndim}\abs{\hat f(y)}(1+\norm{y}^{\ndim+\expb})<+\infty.
\end{eqnarray*}

Of course $\cbsspace\expa\expb\ndim\sset\lspace\ndim$ and
$\sspace\ndim$ is the intersection of all the spaces $\cbsspace\expa\expb\ndim$
when $\expa$ and $\expb$ varies on all the positive real numbers.

Since
\begin{eqnarray*}
	\sspace\ndim&=&\bigcap_{\expa,\expb>0}\cbsspace\expa\expb\ndim
\end{eqnarray*}
is useful to set
\begin{eqnarray*}
	\cbsspace{\infty}{\infty}\ndim&=&\sspace\ndim.
\end{eqnarray*}

Let $\expa,\expb>0$, finite or infinite, be fixed .

Let $f\in\cbsspace\expa\expb\ndim$ be an arbitrary function.

By the very elementary approach to the (Riemann) integration theory
we have
\begin{eqnarray*}\label{eq::riemanint}
	\hat{f}(0)=\int_{\mathbb{R}^n}f(x)\,dx&=&
	\lim_{t\to0^+}t^\ndim\sum_{\ivar\in\Z^\ndim}f(t\ivar).
\end{eqnarray*}

It is then natural to define for $t>0$ and $d>0$
\begin{eqnarray*}
	\cbt{d}{f}{t}&=&\sum_{\ivar\in\Z^\ndim}f(t^{1/d}\ivar),\\
	\cbts{d}{f}{t}&=&{\sum_{\ivar\in\Z^\ndim}}^{'}f(t^{1/d}\ivar)=\cbt{d}{f}{t}-f(0),
\end{eqnarray*}
(where $\sum_{\ivar\in\Z^\ndim}^{'}$,
as usual, stands for $\sum_{\ivar\in\Z^\ndim\setminus\{0\}}$)
and considering the \emph{Mellin transform}.
\begin{eqnarray*}
	\cbx{d}{f}{s}&=&\int_0^{+\infty}\cbts{d}{f}{t}t^s\,\frac{dt}{t}
\end{eqnarray*}
(The introduction of the constant $d$ will simplify some computations).

Our hope is to study the Riemann approximation (\ref{eq::riemanint})
looking at the residues of $\cbx{d}{f}{s}$.

Indeed we have

\begin{theorem}\label{stm::xiA}
Let $\expa,\expb>0$ be given and
let $f\in\cbsspace\expa\expb\ndim$ be an arbitrary function.
Then the integral defining $\cbx{d}{f}{s}$ converges in the strip
$\ndim/d<\Re(s)<(\ndim+\expa)/d$ and the function $\cbx{d}{f}{s}$ extends
to a meromorphic function on the strip $-\expb/d<\Re(s)<(\ndim+\expa)/d$.
having exactly two simple poles respectively at $s=\ndim/d$ with residue $\hat f(0)$
and at $s=0$ with residue $-f(0)$.
Moreover we have the functional equation
\begin{eqnarray*}
	\cbx{d}{f}{s}&=&\cbx{d}{\hat f}{\ndim/d-s},
		\quad-\frac{\expb}{d}<\Re(s)<\frac{\ndim+\expa}{d}.
\end{eqnarray*}
\end{theorem}

The \emph{proof} of the theorem above follows closely the
lines of (some of) the standard proof of the functional equation
for the classical zeta functions used in number theory,
but the \emph{statement} of the theorem in the form above
is not so common.
We refer \cite{article:VSArxivZeta} for a detailed proof.

Observe that $\cbts{d}{f}{t}=\cbts{1}{f}{t^{1/d}}$
and hence, by a simple change of integration variable,
$\cbx{d}{f}{s}=d\cbx{1}{f}{ds}$.
It follows that the function $\cbx{d}{f}{s}$ has a pole at $s=\ndim/d$ if,
and only if, the function $\cbx{1}{f}{s}$ has a pole at $s=\ndim$
with the same residue.

\section{\label{section:ZetaFunctions}Zeta functions machinery}
The following proposition describe the hearth
of our approach to the treatment of Cimmino integrals.

\begin{proposition}\label{stm::RecCompute}
Let $\expa,\expb, d>0$ be given and
let $f\in\sspace\ndim$ be an arbitrary function.
Assume that in the strip $-\expb/d<\Re(s)<(\ndim+\expa)/d$
the function $\cbx{d}{f}{s}$
decomposes as
\begin{eqnarray*}
	\cbx{d}{f}{s}&=&\cbGff(s)\cbZff(s)
\end{eqnarray*}
with $\cbZff(s)$ being meromorphic
with exactly a single simple pole at $s=\ndim/d$.

Then
\begin{eqnarray*}
	\cbZff(0)\Res_{s=0}\cbGff(s)&=&-f(0),\\
	\cbGff(n/d)\Res_{s=\ndim/d}\cbZff(s)&=&\hat f(0).
\end{eqnarray*}
\end{proposition}

\proof
By the previous theorem $\cbx{d}{f}{s}$ has a simple pole at $s=0$
and
\begin{eqnarray*}
	-f(0)&=&\Res_{s=0}\cbx{d}{f}{s}=\Res_{s=0}\cbGff(s)\cbZff(s).
\end{eqnarray*}
Since $\cbZff(s)$ is holomorphic at $s=0$ then
\begin{eqnarray*}
	-f(0)&=&\cbZff(0)\Res_{s=0}\cbGff(s).
\end{eqnarray*}
Replacing $f$ with $\hat f$, we have
\begin{eqnarray*}
	-\hat f(0)&=&\Res_{s=0}\cbx{d}{\hat f}{s}.
\end{eqnarray*}
Using the functional equation $\cbx{d}{f}{s}=\cbx{d}{\hat f}{\ndim/d-s}$ we obtain
\begin{eqnarray*}
	\hat f(0)&=&-\Res_{s=0}\cbx{d}{\hat f}{s}
		=-\Res_{s=0}\cbx{d}{f}{\ndim/d-s}
		=\Res_{s=\ndim/d}\cbx{d}{f}{s}.
\end{eqnarray*}
Since $\cbx{d}{f}{s}$ and $\cbZff(s)$ have a simple pole at $s=\ndim/d$
then necessarily $\cbGff(s)$ is holomorphic at $s=n/d$ and hence
\begin{eqnarray*}
	\hat f(0)&=&\cbGff(n/d)\Res_{s=n/d}\cbZff(s),
\end{eqnarray*}
as desired.

\qed

In the sequel we will apply the above proposition
to functions $f$ which gives
a decompositions $\cbx{d}{f}{s}=\cbGff(s)\cbZff(s)$ where
the function $\cbGff(s)$ is an algebraic combination
of elementary functions and functions of the form
$\Gamma(as+b)$, being $\Gamma(z)$ the Euler gamma function,
and $\cbZff(s)$ is representable when $\Re(s)>>0$ as (generalized)
Dirichlet series
\begin{eqnarray*}
	\cbZff(s)&=&\sum_{k=1}^\infty c_k\lambda_k^{-s},
\end{eqnarray*}
where $c_k,\lambda_k\in\R$,
$\lambda_1<\lambda_2<\cdots$,
and $\lambda_k\to+\infty$ as $k\to+\infty$.

\section{\label{section:QuadraticForms}Quadratic forms}
We begin considering the well known classical example of the
Gaussian integrals associated to quadratic form.

Let us recall that $\Gamma(s)$
denotes the Euler gamma function;
it is a meromorphic function on $C$ 
which for $\Re(s)>0$ satisfies
\begin{eqnarray*}
	\Gamma(s)=\int_0^{+\infty}e^t t^s\frac{dt}{t}.
\end{eqnarray*}
We assume the reader knows all the properties of such a function.

Let $\qA$ be a real symmetric matrix of order $n$.
We set
\begin{eqnarray*}
	\qform{\qA}(x)=\eucdot[\qA x,x]
\end{eqnarray*}

Assume that $\qA$ is positive definite and consider the function
\begin{eqnarray*}
	g_\qA(x)=e^{-\pi\qform{\qA}(x)}
\end{eqnarray*}

then $g_\qA\in\sspace\ndim$ and a standard argument yields
\begin{eqnarray*}
	\hat g_\qA(y)=\frac{1}{\sqrt{\det\qA}}g_{\qA^\negone}(y).
\end{eqnarray*}
In particular
\begin{eqnarray*}
	\hat g_\qA(0)=\frac{1}{\sqrt{\det\qA}}.
\end{eqnarray*}

Since $\qform{\qA}(tx)=t^2\qform{\qA}(x)$ then, choosing $d=2$,
\begin{eqnarray*}
	\cbts{2}{g_\qA}{t}=\cbtsumnz{\ivar}\ndim e^{-t\pi\qform{\qA}(\ivar)}
\end{eqnarray*}
and hence, when $\Re(s)>\ndim/2$,
\begin{eqnarray*}
	\cbx{2}{\gbella{\qA}}{s}&=&\cbtsumnz{\ivar}\ndim\int_0^{+\infty}
		e^{-t\pi\qform{\qA}(\ivar)}t^s\frac{dt}{t}.
\end{eqnarray*}
The change of variable $u=t\pi\qform{\qA}(\ivar)$ yields
\begin{eqnarray*}
	\cbx{2}{\gbella{\qA}}{s}&=&\pi^{-s}\left(\int_0^{+\infty}
		e^{-u}u^s\frac{du}{u}\right)\cbtsumnz{\ivar}\ndim\qform{\qA}(\ivar)^{-s}\\
		&=&\pi^{-s}\Gamma(s)\cbtsumnz{\ivar}\ndim\qform{\qA}(\ivar)^{-s}.
\end{eqnarray*}

The function defined for $\Re(s)>\ndim/2$ by
\begin{eqnarray*}
		\cbza\qA{s}&=&\cbtsumnz{\ivar}\ndim\qform{\qA}(\ivar)^{-s}
\end{eqnarray*}
is the \emph{Epstein zeta function} associated to the
quadratic form $\qform{\qA}$.

Since $\gbella{\qA}(x)\in\sspace\ndim$ then, by Theorem \ref{stm::xiA},
$\cbx{2}{\gbella{\qA}}{s}$ is a meromorphic function on $\C$
having exactly two simple poles at $s=\ndim/2$ and $s=0$ with residues
respectively $R=1$ and $R=-1$.

Since the function $\pi^{-s}\Gamma(s)$ never vanishes on $\C$
and the function $\Gamma(s)$ has a simple pole at $s=0$ with residue $R=1$
it follows that $\cbza\qA{s}$ extends to a meromorphic function on $\C$
having exactly a simple pole at $s=\ndim/2$.

By Proposition \ref{stm::RecCompute} and the relation
\begin{eqnarray*}
	\Gamma(s+1)=s\Gamma(s)
\end{eqnarray*}
we obtain that
the residue of $\cbza\qA{s}$ at $s=\ndim/2$ is
\begin{eqnarray*}
		&&\Res_{s=\ndim/2}\cbza\qA{s}
			=\frac{\pi^{\ndim/2}}{\Gamma\left(\frac\ndim{2}\right)}
			\frac{1}{\sqrt{\det\qA}}
			=\frac{n}{2}
				\left(\frac{\pi^{\ndim/2}}{\Gamma\left(\frac\ndim{2}+1\right)}\right)
				\frac{1}{\sqrt{\det\qA}}
\end{eqnarray*}
and that
\begin{eqnarray*}
		\cbza\qA{s}&=&-1.
\end{eqnarray*}

All that is well known to number theorists.

But now consider two quadratic forms, $\qform{\qA}$ and $\qform{\qB}$ with
$\qA$ and $\qB$ symmetric matrices with $\qA$ positive definite
and set
\begin{eqnarray*}
		\gbellab\qA\qB(x)&=&\qform\qB(x)e^{-\pi\qform\qA(x)}.
\end{eqnarray*}
Then $\gbellab\qA\qB\in\sspace\ndim$.

\begin{proposition}\label{stm::zetaquad}
The Fourier transform of $\gbellab\qA\qB(x)$ is
\begin{eqnarray*}
		\hat\gbell_{\qA,\qB}(y)&=&
		-\frac{1}{\sqrt{\det\qA}}\gbellab{\qA^\negone}{\qA^\negone\qB\qA^\negone}(y)
		+\frac{\Trace(\qA^\negone\qB)}{2\pi\sqrt{\det\qA}}\gbella{\qA^\negone}(y)
		\\
		&=&-\frac{1}{\sqrt{\det\qA}}
			\qform{\qA^\negone\qB\qA^\negone}(y)e^{-\pi\qform{\qA^\negone}(y)}
		+\frac{\Trace(\qA^\negone\qB)}{2\pi\sqrt{\det\qA}}e^{-\pi\qform{\qA^\negone}(y)},
\end{eqnarray*}
where $\Trace(\qA^\negone\qB)$ denotes the trace of the matrix $\qA^\negone\qB$.

In particular we have
\begin{eqnarray*}
		\int_{\R^n}\gbellab\qA\qB(x)\,dx=\hat\gbell_{\qA,\qB}(0)=
			\frac{\Trace(\qA^\negone\qB)}{2\pi\sqrt{\det\qA}}.
\end{eqnarray*}
\end{proposition}

{
\def\inva#1{a_{#1}^{'}}
\def\Asqrt{\qA^{\frac{1}{2}}}
\def\Asqrti{\qA^{-\frac{1}{2}}}
\proof
The function
\begin{eqnarray*}
	&&e^{-\pi\norm{x}^2},
\end{eqnarray*}
where
\begin{eqnarray*}
	&&\norm{x}^2=\eucdot[x,x]=\sum_{i=1}^nx_i^2,
\end{eqnarray*}
coincide with its Fourier transform.
Set
\begin{eqnarray*}
	&&h_{ij}(x)=x_ix_je^{-\pi\norm{x}^2}.
\end{eqnarray*}
Then, the standard properties of the Fourier transform yields
\begin{eqnarray*}
	\hat h_{ij}(y)=-\frac{1}{4\pi^2}D_{y_i}D_{y_j}e^{-\pi\norm{y}^2},
\end{eqnarray*}
where $D_{y_i}$ denote the operator of derivation with respect to
the variable $y_i$.

We compute
\begin{eqnarray*}
	h_{ij}(y)&=&y_iy_je^{-\pi\norm{y}^2},\\
	D_{y_j}h_{ij}(y)&=&-2\pi y_je^{-\pi\norm{y}^2},\\
	D_{y_i}D_{y_j}h_{ij}(y)&=&-2\pi D_{y_i}\bigl(y_je^{-\pi\norm{y}^2}\bigr)
	=4\pi^2e^{-\pi\norm{y}^2}-2\pi \delta_{ij}e^{-\pi\norm{y}^2},\\
	\hat h_{ij}(y)&=&-\frac{1}{4\pi^2}D_{y_i}D_{y_j}e^{-\pi\norm{y}^2}
	=-e^{-\pi\norm{y}^2}+\frac{1}{2\pi}\delta_{ij}e^{-\pi\norm{y}^2}.
\end{eqnarray*}

Let denote by $I_n$ the identity matrix of order $n$.
Let $\qC=(c_{ij})$ be a real symmetric matrix of order $n$ and set
\begin{eqnarray*}
	&&h_\qC(x)=\gbell_{I_n,\qC}(x)
		=\qform\qC(x)e^{-\pi\norm{x}^2}
		=\sum_{i,j=1}^nc_{ij}h_{ij}(x)
\end{eqnarray*}
Then
\begin{eqnarray*}
	\hat h_\qC(y)&=&\sum_{i,j=1}^nc_{ij}\hat h_{ij}(y)\\
	&=&\sum_{i,j=1}^n-c_{ij}e^{-\pi\norm{y}^2}
		+\sum_{i,j=1}^n\frac{1}{2\pi}\qc_{ij}\delta_{ij}e^{-\pi\norm{y}^2}\\
	&=&-\qform\qC(y)e^{-\pi\norm{y}^2}
		+\frac{1}{2\pi}\Trace(\qC)e^{-\pi\norm{y}^2}\\
\end{eqnarray*}
Let $\Asqrt$ be the unique positive definite symmetric matrix
suc that $$\left(\Asqrt\right)^2=\qA$$ and let $\Asqrti$ be its inverse.

Set
\begin{eqnarray*}
	&&f(x)=h_\qC(x)
\end{eqnarray*}
where $\qC=\Asqrti\qB\Asqrti$.
We have
\begin{eqnarray*}
	&&\Trace(\qC)=\Trace(\Asqrti\qB\Asqrti)
		=\Trace(\Asqrti\Asqrti\qB)
		=\Trace(\qA^\negone\qB)
\end{eqnarray*}
and hence
\begin{eqnarray*}
	\hat f(y)=-\qform{\qC}(y)e^{-\pi\norm{y}^2}
		+\frac{1}{2\pi}\Trace(\qA^\negone\qB)e^{-\pi\norm{y}^2}.
\end{eqnarray*}
We also have
\begin{eqnarray*}
	\qform{\qC}(\Asqrt x)&=&\qform\qB(x),\\
	\qform{\qC}(\Asqrti y)&=&\qform{\qA^\negone\qB\qA^\negone}(y),\\
	\norm{\Asqrt y}&=&\qform\qA(y).
\end{eqnarray*}
It follows that
\begin{eqnarray*}
	&\gbellab\qA\qB(x)=f(\Asqrt x)
\end{eqnarray*}
and hence
\begin{eqnarray*}
	\hat\gbell_{\qA,\qB}(y)&=&\frac{1}{\sqrt{\det\qA}}\hat f(\Asqrti y)\\
		&=&-\frac{1}{\sqrt{\det\qA}}
			\qform{\qA^\negone\qB\qA^\negone}(y)e^{-\pi\qform{\qA^\negone}(y)}
		+\frac{\Trace(\qA^\negone\qB)}{2\pi\sqrt{\det\qA}}e^{-\pi\qform{\qA^\negone}(y)},
\end{eqnarray*}
as desired.

\qed

Let now compute $\cbts{2}{\gbellab\qA\qB}{t}$
and $\cbx2{\gbellab\qA\qB}{s}$.

Since $\qform{\qA}(tx)=t^2\qform{\qA}(x)$ and $\qform{\qB}(tx)=t^2\qform{\qB}(x)$
we have
\begin{eqnarray*}
	&&\cbts{2}{\gbellab\qA\qB}{t}=
		\cbtsumnz{\ivar}\ndim t\qform{\qB}(\ivar)e^{-t\pi\qform{\qB}(\ivar)},
\end{eqnarray*}
and
\begin{eqnarray*}
	\cbx{2}{\gbellab\qA\qB}{s}&=&\int_0^{+\infty}
		\cbtsumnz\ivar\ndim t\qform{\qB}(\ivar)e^{-t\pi\qform{\qB}(\ivar)}t^s
		\,\frac{dt}{t}\\
	&=&\cbtsumnz\ivar\ndim\int_0^{+\infty}
		\qform{\qB}(\ivar)e^{-t\pi\qform{\qB}(\ivar)}t^{s+1}
		\,\frac{dt}{t}.
\end{eqnarray*}
By the change of variable $u=t\pi\qform{\qB}(\ivar)$ we obtain
\begin{eqnarray*}
	\cbx{2}{\gbellab\qA\qB}{s}
	&=&\pi^{-(s+1)}\left(\int_0^{+\infty}
		e^{-u}u^{s+1}
		\,\frac{dt}{t}\right)
		\cbtsumnz\ivar\ndim\qform{\qB}(\ivar)\qform{\qA}(\ivar)^{-(s+1)}\\
	&=&\pi^{-(s+1)}\Gamma(s+1)
		\cbtsumnz\ivar\ndim\qform{\qB}(\ivar)\qform{\qA}(\ivar)^{-(s+1)}.
\end{eqnarray*}

Set
\begin{eqnarray*}
	\cbzab\qA\qB{s}=
		\cbtsumnz\ivar\ndim\qform{\qB}(\ivar)\qform{\qA}(\ivar)^{-s}
\end{eqnarray*}

The meromorphic function $\pi^{-(s+1)}\Gamma(s+1)$ never vanishes on $\C$.
Since $\gbellab\qA\qB(0)=0$ by Theorem \ref{stm::xiA} the function
$\cbx{2}{\gbellab\qA\qB}{s}$ is a meromorphic function on $\C$
having exactly a simple pole at $s=\ndim/2$ and hence
the sum defining $\cbzab\qA\qB{s}$
converges for $\Re(s)>\ndim/2$ and extends to a meromorphic
function on $\C$ having exactly a simple pole at $s=\ndim/2$.
Then proposition \ref{stm::RecCompute} and
proposition \ref{stm::zetaquad} give:

\begin{theorem}\label{stm::reszetaquad}
Let $\qA$ and $\qB$ be two symmetric matrices of order $n$
with $\qA$ positive definite. Then the sum defining
the zeta function $\cbzab\qA\qB{s}$
converges absolutely for $\Re(s)>\ndim/2+1$.

The zeta function $\cbzab\qA\qB{s}$
extends to a meromorphic function on $\C$ having
exactly a simple pole at $s=\ndim/2+1$ with residue
\begin{eqnarray*}
	\Res_{s=\ndim/2+1}\cbzab\qA\qB{s}
	=\frac{1}{2}\left(\frac{\pi^{\ndim/2}}{\Gamma\left(\frac\ndim{2}+1\right)}\right)
	\frac{\Trace(\qA^\negone\qB)}{\sqrt{\det\qA}}.
\end{eqnarray*}
\end{theorem}

From the equality
\begin{eqnarray*}
		&&\hat\gbell_{\qA,\qB}(y)=
		-\frac{1}{\sqrt{\det\qA}}\gbellab{\qA^\negone}{\qA^\negone\qB\qA^\negone}(y)
		+\frac{\Trace(\qA^\negone\qB)}{2\pi\sqrt{\det\qA}}\gbella{\qA^\negone}(y)
\end{eqnarray*}
proved in Proposition \ref{stm::zetaquad} we obtain
\begin{eqnarray*}
	\cbx{2}{\hat\gbell_{\qA,\qB}}{s}&=&
		-\frac{1}{\sqrt{\det\qA}}\pi^{-(s+1)}\Gamma(s+1)
			\cbza{\qform{\qA^\negone}}{\qform{\qA^\negone\qB\qA^\negone}}{s+1}\\
		&&\quad+\frac{\Trace(\qA^\negone\qB)}{2\sqrt{\det\qA}}\pi^{-(s+1)}\Gamma(s)
			\cbza{\qform{\qA^\negone}}{s}.
\end{eqnarray*}
Inserting such expressions in the functional equation
\begin{eqnarray*}
	\cbx{2}{\gbellab\qA\qB}{\frac\ndim{2}-s}
	=\cbx{2}{\hat\gbell_{\qA,\qB}}{s}
\end{eqnarray*}
and dividing by $\pi$ we obtain:

\begin{theorem}\label{stm::funeqzetaquad}
Let $\qA$ and $\qB$ be as in theorem \ref{stm::reszetaquad}.
The zeta function $\cbzab\qA\qB{s}$
satisfies the functional equation
\begin{eqnarray*}
	&&\pi^{-\left(\frac\ndim{2}-s\right)}
		\Gamma\left(\frac\ndim{2}+1-s\right)
		\cbzab\qA\qB{\frac\ndim{2}+1-s}\\
	&&\quad+\frac{1}{\sqrt{\det\qA}}\pi^{-s}
		\Gamma\left(s+1\right)
		\cbzab{\qA^\negone}{\qA^\negone\qB\qA^\negone}{s+1}\\
	&&\quad\quad=\frac{\Trace(\qA^\negone\qB)}{2\sqrt{\det\qA}}\pi^{-s}
		\Gamma(s)\cbza{\qA^\negone}{s}.
\end{eqnarray*}
\end{theorem}

%

\section{\label{section:Lattices}Lattices}
{
\def\lU{U}
\def\lAbis{\lA_1}

A \emph{lattice} in $\R^\ndim$ is a set of the form
\begin{eqnarray*}
	\Lattice=\bigl\{\lA\ivar\mid\ivar\in\Z^\ndim\bigr\}
\end{eqnarray*}
where $\lA\in GL(\ndim,\R)$ is a real invertible matrix of order $\ndim$.

If $\lAbis$ is an other matrix such that $\lAbis(\Z^\ndim)=\Lattice$
then $\lU=\lA\lAbis^\negone(\Z^\ndim)\sset\Z^\ndim$.
It follows that the invertible matrix $\lU$ has integral coefficients
and hence $\abs{\det\lU}=1$, that is
$\abs{\det\lA}=\abs{\det\lAbis}$.
We define the volume of the lattice $\Lattice$ as
\begin{eqnarray*}
	\abs\Lattice=\abs{\det\lA}.
\end{eqnarray*}
By the argument given above the definition
does not depend on the choice of the matrix $\lA$.

If $\Lattice=\lA(\Z^\ndim)$ is a lattice the
\emph{dual lattice} is the lattice $\cbdual\Lattice$
associate to the inverse of the transpose of the matrix $A$.
For convenience we also set
\begin{eqnarray*}
\cbmdual\lA=(\cbTranspose\lA)^\negone
\end{eqnarray*}

Given a lattice $\Lattice\sset\R^\ndim$
and a positive definite symmetric matrix $\qA$ 
of order $\ndim$ we define
\begin{eqnarray*}
	\cblza\Lattice\qA{s}=\cblsumnz\ivar\Lattice\qform\qA(\ivar)^{-s}
\end{eqnarray*}
and if $\qB$ is any symmetric matrix we also define
\begin{eqnarray*}
	\cblzab\Lattice\qA\qB{s}=\cblsumnz\ivar\Lattice\qform\qB(\ivar)\qform\qA(\ivar)^{-s}.
\end{eqnarray*}

Of course we have
\begin{eqnarray*}
	\cblza{\Z^\ndim}\qA{s}=\cbza\qA{s}
\end{eqnarray*}
and
\begin{eqnarray*}
	\cblzab{\Z^\ndim}\qA\qB{s}=\cbzab\qA\qB{s}.
\end{eqnarray*}

If $\Lattice=\lA(\Z^\ndim)$ with $\lA\in GL(\ndim,\R)$ then we
have
\begin{eqnarray*}
	\cblza\Lattice\qA{s}=\cbza{\cbTranspose\lA\qA\lA}{s},
\end{eqnarray*}
where $\cbTranspose\lA$ denotes the transpose of the matrix $\lA$,
and
\begin{eqnarray*}
	\cblzab\Lattice\qA\qB{s}=
		\cbzab{\cbTranspose\lA\qA\lA}{\cbTranspose\lA\qB\lA}{s}.
\end{eqnarray*}
We also have
\begin{eqnarray*}
	&&\sqrt{\det(\cbTranspose\lA\qA\lA)}=\abs{\det A}\sqrt{\det\qA}=
		\abs\Lattice\sqrt{\det A}
\end{eqnarray*}
and
\begin{eqnarray*}
	\Trace\bigr((\cbTranspose\lA\qA\lA)^\negone(\cbTranspose\lA\qB\lA)
		(\cbTranspose\lA\qA\lA)^\negone\bigl)
		&=&\Trace\bigr(\lA^\negone(\qA^\negone\qB\qA^\negone)\lA\bigr)\\
		&=&\Trace(\qA^\negone\qB\qA^\negone).
\end{eqnarray*}
Applying the theta-zeta machinery to the function
\begin{eqnarray*}
		\lgbellab\lA\qA\qB(x)
		=\gbellab\qA\qB(\lA x)&=&\qform\qB(\lA x)e^{-\pi\qform\qA(\lA x)},
\end{eqnarray*}
observing that the Fourier transform of $\lgbellab\lA\qA\qB(x)$
is
\begin{eqnarray*}
		\hat\gbell_{\lA,\qA,\qB}(y)=
		\frac{1}{\abs{\det\lA}}\hat\gbell_{\qA,\qB}(\cbmdual\lA y),
\end{eqnarray*}
it follows that the results of the previous section generalize:

\begin{theorem}\label{stm::reslzetaquad}
Let $\Lattice\sset\R^\ndim$ be a lattice and
Let $\qA$ and $\qB$ be two symmetric matrices of order $n$
with $\qA$ positive definite. Then the sum defining
the zeta functions $\cblza\Lattice\qA{s}$ and
$\cblzab\Lattice\qA\qB{s}$
converges absolutely respectively for $\Re(s)>\ndim/2$
and $\Re(s)>\ndim/2+1$.

The zeta functions $\cblza\Lattice\qA{s}$ and $\cblzab\Lattice\qA\qB{s}$
extends to a meromorphic function on $\C$ having
exactly a simple pole respectively at $s=\ndim/2$ and $s=\ndim/2+1$ 
with residues
\begin{eqnarray*}
	&&\Res_{s=\ndim/2}\cblza\Lattice\qA{s}
		=\frac{n}{2}
			\left(\frac{\pi^{\ndim/2}}{\Gamma\left(\frac\ndim{2}+1\right)}\right)
			\frac{1}{\abs\Lattice\sqrt{\det\qA}},
\end{eqnarray*}
and			
\begin{eqnarray*}
	&&\Res_{s=\ndim/2+1}\cblzab\Lattice\qA\qB{s}
		=\frac{1}{2}\left(\frac{\pi^{\ndim/2}}
			{\Gamma\left(\frac\ndim{2}+1\right)}\right)
			\frac{\Trace(\qA^\negone\qB)}{\abs\Lattice\sqrt{\det\qA}}.
\end{eqnarray*}
\end{theorem}

\begin{theorem}\label{stm::funeqlzetaquad}
Let $\Lattice\sset\R^\ndim$, $\qA$ and $\qB$ be as in theorem \ref{stm::reslzetaquad}.
The zeta functions $\cblza\Lattice\qA{s}$ and $\cblzab\Lattice\qA\qB{s}$
satisfy the functional equations respectively
\begin{eqnarray*}
	&&\pi^{-\left(\frac\ndim{2}-s\right)}
		\Gamma\left(\frac\ndim{2}-s\right)
		\cblza\Lattice\qA{\frac\ndim{2}-s}\\
	&&\quad=\frac{1}{\abs\Lattice\sqrt{\det\qA}}\pi^{-s}
		\Gamma(s)\cblza{\cbdual\Lattice}{\qA^\negone}{s}.
\end{eqnarray*}
and
\begin{eqnarray*}
	&&\pi^{-\left(\frac\ndim{2}-s\right)}
		\Gamma\left(\frac\ndim{2}+1-s\right)
		\cblzab\Lattice\qA\qB{\frac\ndim{2}+1-s}\\
	&&\quad+\frac{1}{\abs\Lattice\sqrt{\det\qA}}\pi^{-s}
		\Gamma\left(s+1\right)
		\cblzab{\cbdual\Lattice}{\qA^\negone}{\qA^\negone\qB\qA^\negone}{s+1}\\
	&&\quad\quad=\frac{\Trace(\qA^\negone\qB)}{2\abs\Lattice\sqrt{\det\qA}}\pi^{-s}
		\Gamma(s)\cblza{\cbdual\Lattice}{\qA^\negone}{s}.
\end{eqnarray*}

\end{theorem}

}

%

\section{\label{section:Integrals}Integral representation}

Let $\qform{\qA}$ and $\qform{\qB}$ be two quadratic form with
$\qA$ and $\qB$ symmetric matrices and $\qA$ positive definite.

We now will give such residues as integrals over
the boundary of the unit ball in $\R^\ndim$.

We already observed that 
\begin{eqnarray*}
	\int_{\R^\ndim}e^{-\pi\qform\qA(x)}\,dx=\frac{1}{\sqrt{\det\qA}}
\end{eqnarray*}
and by Proposition \ref{stm::zetaquad}
we also have
\begin{eqnarray*}
	\int_{\R^\ndim}\qform\qB(x)e^{-\pi\qform\qA(x)}\,dx=
		\frac{\Trace(\qA^\negone\qB)}{2\pi\sqrt{\det\qA}}.
\end{eqnarray*}

If $\lA\in GL(\ndim,\R)$ then a simple change of
variable gives
\begin{eqnarray*}
	\int_{\R^\ndim}e^{-\pi\qform\qA(\lA x)}\,dx
		=\frac{1}{\abs{\det\lA}\sqrt{\det\qA}}
\end{eqnarray*}
and
\begin{eqnarray*}
	\int_{\R^\ndim}\qform\qB(\lA x)e^{-\pi\qform\qA(\lA x)}\,dx=
		\frac{\Trace(\qA^\negone\qB)}{2\pi\abs{\det\lA}\sqrt{\det\qA}}.
\end{eqnarray*}

Let us recall that given $f\in L^1(\R^\ndim)$ the
integration by polar coordinates gives
\begin{eqnarray*}
	\int_{\R^\ndim}f(x)\,dx=
	\int_{S^{\ndim-1}}\left(
		\int_0^{+\infty}f(ru)r^{n-1}\,dr
	\right)\,du,
\end{eqnarray*}
where $du$ is the Euclidean (hyper-)surface measure on
the unit sphere $S^{\ndim-1}$.

Using such formula we obtain
\begin{eqnarray*}
	\frac{\Trace(\qA^\negone\qB)}{2\pi\abs{\det\lA}\sqrt{\det\qA}}
	&=&\int_{\R^\ndim}\qform\qB(\lA x)e^{-\pi\qform\qA(\lA x)}\,dx\\
	&=&
		\int_{S^{\ndim-1}}\left(
			\int_0^{+\infty}\qform\qB(r\lA u)e^{-\pi\qform\qA(r\lA u)}r^{n-1}\,dr
		\right)\,du\\
	&=&
		\int_{S^{\ndim-1}}\qform\qB(\lA u)\left(
			\int_0^{+\infty}e^{-\pi\qform\qA(\lA u)}r^{n+2}\frac{dr}{r}
		\right)\,du
\end{eqnarray*}
By the substitution $\pi r^2\qform\qA(\lA u)=t$ in the inner integral
we obtain
\begin{eqnarray*}
	&&\int_{S^{\ndim-1}}\qform\qB(\lA u)\left(
		\int_0^{+\infty}e^{-\pi\qform\qA(\lA u)}r^{n+2}\frac{dr}{r}
	\right)\,du\\
	&&=
		\frac{1}{2}
		\int_{S^{\ndim-1}}\qform\qB(\lA u)\left(
			\int_0^{+\infty}e^{-t}
				\left(\frac{t}{\pi\qform\qA(\lA u)}\right)^{n/2+1}
				\frac{dt}{t}
		\right)\,du\\
	&&=
		\frac{1}{2}\pi^{-(n/2+1)}
		\Gamma\left(\frac{n}{2}+1\right)
		\int_{S^{\ndim-1}}\qform\qB(\lA u)\qform\qA(\lA u)^{-(n/2+1)}\,du
\end{eqnarray*}
and hence
\begin{eqnarray*}		
		\int_{S^{\ndim-1}}\qform\qB(\lA u)\qform\qA(\lA u)^{-(n/2+1)}\,du=
		\left(\frac{\pi^{n/2}}{\Gamma\left(\frac{n}{2}+1\right)}\right)
		\frac{\Trace(\qA^\negone\qB)}{\abs{\det\lA}\sqrt{\det\qA}}.
\end{eqnarray*}

When $\qB=\qA$ we obtain
\begin{eqnarray*}		
		\int_{S^{\ndim-1}}\qform\qA(\lA u)^{-n/2}\,du=
		\ndim\left(\frac{\pi^{n/2}}{\Gamma\left(\frac{n}{2}+1\right)}\right)
		\frac{1}{\abs{\det\lA}\sqrt{\det\qA}}.
\end{eqnarray*}

Replacing $s$ with $s/2$ in the formulas
of the residues in Theorem \ref{stm::reslzetaquad} we obtain
\begin{eqnarray*}
	&&\Res_{s=\ndim}\cblza\Lattice\qA{s/2}
		=\ndim
			\left(\frac{\pi^{\ndim/2}}{\Gamma\left(\frac\ndim{2}+1\right)}\right)
			\frac{1}{\abs\Lattice\sqrt{\det\qA}},
\end{eqnarray*}
and			
\begin{eqnarray*}
	&&\Res_{s=\ndim+2}\cblzab\Lattice\qA\qB{s/2}
		=\left(\frac{\pi^{\ndim/2}}
			{\Gamma\left(\frac\ndim{2}+1\right)}\right)
			\frac{\Trace(\qA^\negone\qB)}{\abs\Lattice\sqrt{\det\qA}}.
\end{eqnarray*}

Thus we obtained:

\begin{theorem}\label{stm::integralrep}
Let $\qform{\qA}$ and $\qform{\qB}$ be two quadratic form with
$\qA$ and $\qB$ symmetric matrices and $\qA$ positive definite.
Let $\lA\in GL(\ndim,\R)$ and set $\Lattice=\lA(\Z^\ndim)$.
Then we have
\begin{eqnarray*}
	&&\Res_{s=\ndim}\cblza\Lattice\qA{s/2}
		=\int_{S^{\ndim-1}}\qform\qA(\lA u)^{-n/2}\,du.
\end{eqnarray*}
and
\begin{eqnarray*}
	&&\Res_{s=\ndim+2}\cblzab\Lattice\qA\qB{s/2}
		=\int_{S^{\ndim-1}}\qform\qB(\lA u)\qform\qA(\lA u)^{-(n/2+1)}\,du,
\end{eqnarray*}
\end{theorem}

%

\section{\label{section:LinearSystems}Linear systems}
{
We recall that $\cbmdual\lA$ denotes the inverse of the
transpose of the matrix $\lA$.

Given $\lA\in GL(\R^\ndim)$ we define
\begin{eqnarray*}
	\cbcza\lA{s}=\cbtsumnz\ivar\ndim
		\norm{\lA\om}^{-2s}
\end{eqnarray*}
and given also $b\in\R^\ndim$
we define the vector value zeta function
\begin{eqnarray*}
	\cbczab\lA{b}{s}=\cbtsumnz\ivar\ndim
		\norm{\lA\om}^{-2s}\eucdot[b,\om]\lA\om.
\end{eqnarray*}

Of course we have
\begin{eqnarray*}
	\cbcza\lA{s}=\cblza\Lattice{I_\ndim}{s}
\end{eqnarray*}
with $\Lattice=\lA(\Z^\ndim)$, $I_\ndim$
the identity matrix of order $\ndim$
and hence,
by Theorem \ref{stm::reslzetaquad},
we have
\begin{eqnarray*}
		&&\Res_{s=\ndim/2}\cbcza\lA{s}
			=\frac{n}{2}
				\left(\frac{\pi^{\ndim/2}}{\Gamma\left(\frac\ndim{2}+1\right)}\right)
				\frac{1}{\abs{\det\lA}}
\end{eqnarray*}

But we also have:
\begin{theorem}\label{stm::reslinear}
The series defining $\cbczab\lA{b}{s}$ converges for $\Re(s)>\ndim/2+1$
and the function $\cbczab\lA{b}{s}$ extends to a (vector value) 
meromorphic function on $\C$ having only a simple pole at $s=\ndim/2+1$
with residue
\begin{eqnarray*}
	\Res_{s=\ndim/2+1}\cbczab\lA{b}{s}
	=\frac{1}{2}\left(\frac{\pi^{\ndim/2}}{\Gamma\left(\frac\ndim{2}+1\right)}\right)
	\frac{1}{\abs{\det\lA}}\cbmdual\lA b.
\end{eqnarray*}

Moreover, if $c\in\R^\ndim$ then we have the functional equation
\begin{eqnarray*}
	&&\pi^{-\left(\frac\ndim{2}-s\right)}
		\Gamma\left(\frac\ndim{2}+1-s\right)
		\eucdot[\cbczab\lA{b}{\frac\ndim{2}+1-s},\lA c]\\
	&&\quad+\frac{\pi^{-s}}{\abs{\det\lA}}
		\Gamma\left(s+1\right)
		\eucdot[\cbmdual\lA b,\cbczab{\cbmdual\lA}{c}{s+1}]\\
	&&\quad\quad=\frac{\eucdot[b,c]}{2\abs{\det\lA}}\pi^{-s}
		\Gamma(s)\cbcza{\cbmdual\lA}{s}.
\end{eqnarray*}
\end{theorem}

\proof {
\def\avec{u}
\def\bvec{v}
\def\cvec{w}
Given $\avec,\bvec\in\R^\ndim$, with $\avec=(\avec_1,\ldots,\avec_\ndim)$ and $\bvec=(\bvec_1,\ldots,\bvec_\ndim)$
we denote by $\avec\otimes\bvec$ the symmetric matrix $(\cvec_{ij})$
of order $\ndim$ with entries given by $\cvec_{ij}=\avec_i\bvec_j$.
Then we have
\begin{eqnarray*}
	\Trace(\avec\otimes\bvec)=\eucdot[\avec,\bvec].
\end{eqnarray*}

Now observe that
\begin{eqnarray*}
	&&\eucdot[\cbczab\lA{b}{s},\lA c]=
		\cblzab\Lattice{I_n}\qB{s}
\end{eqnarray*}
where $\Lattice=\lA(\Z^n)$, $I_n$ is the identity matrix of order $n$
and
\begin{eqnarray*}
	&&\qB=\cbmdual\lA b\otimes\lA c.
\end{eqnarray*}
We have
\begin{eqnarray*}
	&&\Trace(\cbmdual\lA b\otimes\lA c)
		=\eucdot[\cbmdual\lA b,\lA c]
		=\eucdot[b,c]
\end{eqnarray*}
and hence, the formulas for the residue and
the functional equation of the function 
$\eucdot[\cbczab\lA{b}{\frac\ndim{2}+1-s},\lA c]$
follows from Theorem \ref{stm::reslzetaquad}
and Theorem \ref{stm::funeqlzetaquad}

}\qed

We are now ready to state and prove the reformulation
of Cimmino's results.

\begin{theorem}\label{stm::ResSolve}
Let
\begin{eqnarray*}
\lA x&=&b
\end{eqnarray*}
be a linear system of $\ndim$ equation and $\ndim$ unknown,
where the unknown values $x_1,\ldots,x_\ndim$ are
the components if the column vector $x\in\R^\ndim$
and $\lA$ is non singular matrix with real coefficient of
order $\ndim$.

Then we have
\begin{eqnarray*}
x_i&=&\frac{R_i}{R},\quad i=1,\ldots\ndim,
\end{eqnarray*}
where
\begin{eqnarray*}
	R&=&\Res_{s=\ndim}\cbcza{\cbTranspose\lA}{s/2}
\end{eqnarray*}
and
\begin{eqnarray*}
	R_i&=&\ndim\Res_{s=\ndim+2}\eucdot[\cbczab{\cbTranspose\lA}{b}{s/2},e_i]
	,\quad i=1,\ldots\ndim,
\end{eqnarray*}
where $e_1,\ldots,e_\ndim$ is the canonical basis of $\R^\ndim$.

Moreover we have the identities
\begin{eqnarray*}
	R&=&
		\int_{S^{\ndim-1}}\norm{\cbTranspose\lA u}^{-\ndim}\,du,
\end{eqnarray*}
and
\begin{eqnarray*}
	R_i&=&
		\ndim\int_{S^{\ndim-1}}
			\norm{\cbTranspose\lA u}^{-\ndim-2}
				\eucdot[b,u]\eucdot[\cbTranspose\lA x,e_i]\,du.
\end{eqnarray*}

\end{theorem}

\proof
By assumption
\begin{eqnarray*}
	x=\lA^\negone b
\end{eqnarray*}
and we have
\begin{eqnarray*}
	x_i=\eucdot[x,e_i],\quad i=1,\ldots,\ndim
\end{eqnarray*}
where $e_1,\ldots,e_\ndim$ is the canonical basis of $\R^\ndim$.
As previously observed we have
\begin{eqnarray*}
		&&R
			=\ndim
				\left(\frac{\pi^{\ndim/2}}{\Gamma\left(\frac\ndim{2}+1\right)}\right)
				\frac{1}{\abs{\det{\cbTranspose\lA}}}
			=\ndim
				\left(\frac{\pi^{\ndim/2}}{\Gamma\left(\frac\ndim{2}+1\right)}\right)
				\frac{1}{\abs{\det\lA}}
\end{eqnarray*}
and for $i=1,\ldots,\ndim$ we have
\begin{eqnarray*}
	R_i
	&=&\ndim\frac{1}{2}\left(\frac{\pi^{\ndim/2}}
		{\Gamma\left(\frac\ndim{2}+1\right)}\right)
		\frac{1}{\abs{\det{\cbTranspose\lA}}}
		\eucdot[\cbTranspose{\cbmdual\lA} b,e_i]\\
	&=&\frac{\ndim}{2}\left(\frac{\pi^{\ndim/2}}
		{\Gamma\left(\frac\ndim{2}+1\right)}\right)
		\frac{1}{\abs{\det\lA}}
		\eucdot[\lA^\negone b,e_i]
\end{eqnarray*}
and hence
\begin{eqnarray*}
	\frac{R_i}{R}=\eucdot[\lA^\negone b,e_i]=x_i,
\end{eqnarray*}
as required.

The last assertions of the Theorem follows from
Theorem \ref{stm::integralrep}.

\qed

Our tour around Cimmino's ideas is completed.
}

%

\nocite{article:CimminoA}
\nocite{article:CimminoB}
\nocite{article:CimminoC}
\nocite{article:CimminoD}
\nocite{article:CimminoE}
\nocite{book:CimminoOpereScelte}
\nocite{article:BenziPerCimmino}

\end{document}